\documentclass[11pt, a4paper]{article}
\usepackage[centertags]{amsmath}
\usepackage{amsfonts}
\usepackage{amssymb}
\usepackage{amsthm}
\usepackage{newlfont}

\newlength{\defbaselineskip}
\theoremstyle{plain}
\newtheorem{thm}{Theorem}[section]
\newtheorem{lem}{Lemma}[section]

\numberwithin{equation}{section}

\newcommand{\A}{\alpha}
\newcommand{\B}{\beta}

\newcommand{\M}{\overline{M}}
\newcommand{\MM}{\widetilde{M}}
\newcommand{\g}{\overline{g}}
\newcommand{\G}{\widetilde{g}}

\newcommand{\PB}{\overline{P}}
\newcommand{\PP}{\widetilde{P}}

\title{$(a, 1)f$ structures on product of spheres}

\author{Cristina-Elena Hre\c{t}canu\\
\c{S}tefan cel Mare University, Suceava, Rom\^ania\\
Mircea Crasmareanu\\
Al. I. Cuza University, Ia\c{s}i, Rom\^ania}
\date{}

\usefont{T1}{cmss}{m}{ol}

\begin{document}
\maketitle

\begin{abstract}
Our aim in this paper is to give some examples of $(a,  1)f$
Riemannian structures (a generalization of an $r$-paracontact
structure) induced on product of spheres of codimension $r$ ($r \in
\{1,2\} $) in an $m$-dimensional Euclidean space ($m>2$),
endowed with an almost product structure. \\
\end{abstract}

{\bf Keywords}: differential structures, product of spheres, induced structures on submanifolds \\
{\bf MSC 2000}: 53B25,53C15.

\section*{Introduction}
\normalfont

A classification of smooth structures on product of spheres of the
form $S^{k} \times S^{p}$, where $2 \leq k \leq p$ and $ k+p \geq
6$, was given by R. De Sapio \cite{Sapio}. In \cite{Ajala},
S.O.Ajala extended De Sapio's result to smooth structures on $S^{p}
\times S^{q} \times S^{r}$ where $2\leq p \leq q \leq r$. Also, a
complete classification of smooth structures on a generalized
product of spheres was given in \cite{Ajala2}.

By studying properties of some structures constructed on Riemannian
manifolds \cite{Adati, Anastasiei, Goldberg, Ianus, Pitis, Yano}, we
obtain a generalization of $r$-paracontact structure, constructed as
an induced structure on a submanifold in an almost product
Riemannian manifold.

In this paper we show that, if $M$ is a submanifold of codimension
1, isometrically immersed in $\M$, and $\M$ is also of codimension 1
and isometrically immersed in an n-dimensional
 almost product Riemannian manifold $(\MM,\G,\PP)$ ($n>2$), so
that $(M,g)\hookrightarrow (\M,\g) \hookrightarrow(\MM,\G)$ then,
the induced $(a,1)f$ structure on $M$ by the structure $(\PP,\G)$
from $\MM$ is a structure of type
$(P,g,u_{1},u_{2},\xi_{1},\xi_{2}^{\top}, (a_{\A\B}))$, which is the
same that one induced on $M$ by the structure
$(\overline{P},\g,u_{2},\xi_{2},a_{22})$ induced on $\M$ by the
structure from $\MM$.

Finally we give some examples for induced $(a,1)f$ Riemannian
structures on product of spheres of codimension $r$ ($r \in \{1,2\}
$) in an Euclidean space of dimension $m>2$ endowed with an almost
product structure.

\section{Submanifolds in almost product Riemannian \\ manifolds}

Let $\MM$ be an $m$-dimensional Riemannian manifold endowed with a
pair $(\widetilde{P},\widetilde{g})$ where $\widetilde{g}$ is a
Riemannian metric and $\widetilde{P}$ is an $(1,1)$ tensor field so
that $\widetilde{P}^{2}=\varepsilon Id$ for $\varepsilon \in
\{1,-1\}$. We suppose that $\G$ and $\PP$ verify the compatibility
condition
$\widetilde{g}(\widetilde{P}U,\widetilde{P}V)=\widetilde{g}(U,V)$
for every $U, V \in \chi(\MM)$ where $\chi(\MM)$ is the Lie algebra
of the vector fields on $\MM$. This conditions
 is equivalent with $\G(\PP U, V)=\varepsilon \G(U, \PP V)$
for every $U, V \in \chi(\MM)$.

For $\varepsilon =1$, we obtain that $(\MM,\G,\PP)$ is an {\it
almost product Riemannian manifold}.

Let $M$ be an $n$-dimensional submanifold of codimension $r$ $(n, r
\in \mathbb{N}^{*})$ in an almost product Riemannian manifold
$(\MM,\G,\PP)$ and let $g$ be the Riemannian metric induced on $M$
by $\G$.

If $(N_{1},...,N_{r}):=(N_{\A})$ is a local orthonormal basis in the
normal space of $M$ in $x$, denoted $T_{x}^{\bot}(M)$, for every $x
\in M$ (with $\A \in \{1,...,r\}$) then, decompositions of the
vector fields $\PP X$ and $\PP N_{\A}$, respectively, in the
tangential and normal components on the submanifold $M$ in $\MM$ are
as follows:
\[
 \PP X = P X + \sum_{\A=1}^{r}u_{\A}(X)N_{\A}, \leqno(1.1)
\]
and
\[\PP N_{\A}= \varepsilon \xi_{\A} + \sum_{\B=1}^{r}a_{\A \B} N_{\B}
\leqno(1.2)
\]
for every $X \in \chi(M)$ and $\A \in \{1,...,r\}$.

We called in \cite{HreFr} an $(a, \varepsilon)f$ Riemannian
structure on $M$, induced by $\PP$ from $(\MM,\G)$, the following
data which results from the relations (1.1) and (1.2): $(P, g,
\varepsilon\xi_{\A}, u_{\A}, (a_{\A\B})_{r})$ where $a$ is a
notation for the matrix $(a_{\A\B})_{r}$ and $f:=P$. Here, $P$ is an
$(1, 1)$-tensor field on $M$, $\xi_{\A}$ are tangent vector fields
on $M$, $u_{\A}$ are $1$-forms on $M$ and $(a_{\A\B})_{r}$ is a $r
\times r$ matrix of real functions on $M$. Some properties of
$(a,\varepsilon)f$ Riemannian structures are given by the first
author in \cite{HreFr, Hre}. The $(a,1)f$ Riemannian structure
generalizes the Riemannian almost $r$-paracontact structure
\cite{Bucki} obtained from $(a,1)$f structure for $a=0$, and it was
also considered by T. Adati in \cite{Adati}. A similar structure
induced on $M$ by an almost Hermitian structure on $\MM$ was studied
by K. Yano and M. Okumura \cite{Yano}.

\section{$(a, 1)f$ induced structures on submanifolds \\
in submanifolds of almost product Riemannian manifolds}

In the following statements we suppose that $(\MM,\G,\PP)$ is an
$n$-dimensional ($n>2$) almost product Riemannian manifold and
$(\M,\g)$ is a submanifold of codimension $1$, isometrically
immersed in $\MM$ (with the induced metric $\g$ on $\M$ by $\G$).
Let $N_{1}$ be an unit vector field, normal on $\M$ in $\MM$. Then,
we suppose that $(M,g)$ is a submanifold of codimension $1$
isometrically immersed in $\M$ and let $N_{2}$ be an unit vector
field, normal on $M$ in $\M$. Thus, $(M,g)$ is a submanifold of
codimension $2$ in $(\MM,\G)$ and we have the following isometric
immersions between two Riemannian manifolds: $(M,g)\hookrightarrow
(\M,\g) \hookrightarrow(\MM,\G)$
 and $(N_{1},N_{2})$ is a local orthonormal basis in
$T_{x}^{\bot}(M)$ for every $x \in M$.

From the decompositions in tangential and normal components at
$\M$ in $\MM$ of vector fields $\widetilde{P}\overline{X}$
($\overline{X} \in \chi(\M)$) and $\widetilde{P}N_{1}$
respectively, we obtain:
\[
\widetilde{P}\overline{X}=\overline{P}\:\overline{X}+u_{1}(\overline{X})N_{1},
\quad \widetilde{P}N_{1}=\xi_{1}+a_{11}N_{1}, \leqno(2.1)
\]
for any $\overline{X} \in \chi(\overline{M})$ where $\overline{P}$
is an $(1, 1)$ tensor field on $\M$, $u_{1}$ is an 1-form on $\M$,
$\xi_{1}$ is a tangent vector field on $\M$ and $a_{11}$ is a real
function on $\M$.

\lem  The almost product Riemannian structure $(\PP,\G)$ on a
manifold $\MM$ induces, on any submanifold $\M$ of codimension $1$
in $\MM$, an $(\overline{a}, 1)f$ Riemannian structure, which is a
$(\overline{P}, \g, u_{1}, \xi_{1}, a_{11})$ Riemannian structure,
(with $\overline{a}:=a_{11}$ and $f:=\overline{P}$), where
$\overline{P}$ is an $(1, 1)$ tensor field on $\M$, $u_{1}$ is an
$1$-form on $\M$, $\xi_{1}$ is a tangent vector field on $\M$ and
$a_{11}$ is a real function on $\M$. This structure has the
following properties:
\[
\begin{cases}
(i)\quad
\overline{P}^{2}\overline{X}=\overline{X}-u_{1}(\overline{X})\xi_{1},
\: (\forall) \overline{X} \in \chi(\M),\\
(ii) \quad
u_{1}(\overline{P}\:\overline{X})=-a_{11}u_{1}(\overline{X}),
\:(\forall) \overline{X} \in \chi(\M),\\
(iii)\quad u_{1}(\xi_{1})=1-a_{11}^{2},\\
(iv)\quad \overline{P}\xi_{1}=-a_{11}\xi_{1},\\
\end{cases}
\leqno(2.2)
\]
and
\[
\begin{cases}
(i) \quad u_{1}(\overline{X})=\g(\overline{X},\xi_{1}),
\: (\forall) \overline{X} \in \chi(\M),\\
(ii) \quad
\g(\overline{P}\:\overline{X},\overline{Y})=\g(\overline{X},\overline{P}\:\overline{Y}),
\: (\forall) \overline{X},\overline{Y} \in \chi(\M),\\
(iii)\quad
\g(\overline{P}\:\overline{X},\overline{P}\:\overline{Y})=
\g(\overline{X},\overline{Y}), \: (\forall)
\overline{X},\overline{Y} \in \chi(\M).\quad \square
\end{cases}
\leqno(2.3)
\]

\normalfont  The decompositions in tangential and normal
components on $M$ in $\M$ of the vector fields $\overline{P}X$ ($X
\in \chi(M)$) and $\overline{P}N_{2}$ are, respectively, as
follows:
\[
\overline{P}X=PX+u_{2}(X)N_{2}, \quad
\overline{P}N_{1}=\xi_{2}+a_{22}N_{2}, \leqno(2.4)
\]
for any $X \in \chi(M)$, where $P$ is an $(1, 1)$ tensor field on
$M$, $u_{2}$ is an 1-form on $M$, $\xi_{2}$ is a tangential vector
field on $M$ and $a_{22}$ is a real function on $M$.

On the other hand, we remark that the decomposition of the vector
field $\xi_{1} \in \chi(\M)$ in tangential and normal components
on M in $\M$ has the form $\xi_{1}=\xi_{1}^{\top}+\xi_{1}^{\bot}$
and $\xi_{1}^{\bot}$ and $N_{2}$ are collinear.

\lem  The decompositions in the tangential and normal parts on M
in $\MM$ of vector fields $\widetilde{P}X$ ($X \in \chi(M)$),
$\widetilde{P}N_{1}$ and $\widetilde{P}N_{2}$ are, respectively,
as follows:
\[
\begin{cases}
(i)\quad \widetilde{P}X=PX+u_{1}(X)N_{1}+u_{2}(X)N_{2},\:(\forall) X \in \chi(M)\\
(ii) \quad \widetilde{P}N_{1}=\xi_{1}^{\top}+a_{11}N_{1}+a_{12}N_{2},\\
(iii)\quad \widetilde{P}N_{2}=\xi_{2}+a_{21}N_{1}+a_{22}N_{2},
\end{cases}
\leqno(2.5)
\]
where $P$ is an $(1, 1)$ tensor field on $M$, $u_{1}, u_{2}$ are
$1$-forms on $M$, $\xi_{1}^{\top}, \xi_{2}$ are tangent vector
fields on $M$, $(a_{\A\B})$ (with $\A, \B \in \{1,2\}$) is an $2
\times 2$ matrix, and its entries $a_{11}$, $a_{22}$ and
$a_{12}=a_{21}=\G(\xi_{1}^{\bot},N_{2})$ are real functions on $M$.

\lem The structure $(\overline{P}, \g, \xi_{2}, u_{2}, a_{22})$
(induced on a submanifold $(\M,\g)$ of codimension $1$ in a
$n$-dimensional ($n>2$) almost product Riemannian manifold
$(\MM,\G,\PP)$) also induces, on a submanifold $(M,g)$ of
codimension $1$ in $\M$, a Riemannian structure $(P, g, u_{1},
u_{2}, \xi_{1}, \xi_{2}^{\top}, (a_{\A\B}))$ (where $P, u_{1},
u_{2}, \xi_{1}, \xi_{2}^{\top}$, $(a_{\A\B})$ were defined in the
last two propositions) which has the following properties:
\[
 \begin{cases} (i) \quad
P^{2}X=X-u_{1}(X)\xi_{1}-u_{2}(X)\xi_{2}^{\top},\:(\forall) X \in
\chi(M)\\
(ii)\quad u_{1}(PX)=-a_{11}u_{1}(X)-a_{12}u_{2}(X),\:(\forall) X
\in \chi(M)\\
(iii)\quad u_{2}(PX)=-a_{21}u_{1}(X)-a_{22}u_{2}(X),\:(\forall) X \in \chi(M)\\
(iv)\quad u_{1}(\xi_{1})=1-a_{11}^{2}-a_{12}^{2},\\
(v)\quad u_{2}(\xi_{1})=-a_{11}a_{12}-a_{12}a_{22},\\
(vi)\quad u_{1}(\xi_{2}^{\top})=-a_{11}a_{12}-a_{12}a_{22},\\
(vii)\quad u_{2}(\xi_{2}^{\top})=1-a_{12}^{2}-a_{22}^{2},\\
(viii)\quad P(\xi_{1})=-a_{11}\xi_{1}-a_{12}\xi_{2}^{\top},\\
(ix)\:quad (\xi_{2}^{\top})=-a_{12}\xi_{1}-a_{22}\xi_{2}^{\top},
\end{cases}  \leqno(2.6)
\]
 and the properties which depends on the metric g are:
\[ \begin{cases}
(i)\quad u_{1}(X)=g(X,\xi_{1}),\\
(ii)\quad u_{2}(X)=g(X,\xi_{2}^{\top}),\\
(iii)\quad g(PX,Y)=g(X,PY),\\
(iv)\quad g(PX,PY)=g(X,Y)-u_{1}(X)u_{1}(Y)-u_{2}(X)u_{2}(Y),
\end{cases}
\leqno(2.7)
\]
for any $X,Y \in \chi(M).$

\begin{proof} From $\PP(\PP X)=X$ it
follows that:
\[\PP(PX+u_{1}(X)N_{1}+u_{2}(X)N_{2})=X\],
thus we have:
\[P^{2}X+u_{1}(PX)N_{1}+u_{2}(PX)N_{2}+u_{1}(X)(\xi_{1}+a_{11}N_{1}+a_{12}N_{2})+\]
\[+u_{2}(X)(\xi_{2}^{\top}+a_{12}N_{1}+a_{22}N_{2})=X\] Identifying
the tangential and respectively, normal components on M from the
last equality, we obtain (i), (ii) and (iii) from $(2.6)$.

On the other hand, from $\PP(\PP N_{1})=N_{1}$ we derive:
\[N_{1}=\PP (\PP N_{1})=\PP (\xi_{1}+a_{11}N_{1}+a_{12}N_{2})=\]
\[=P\xi_{1}+u_{1}(\xi_{1})N_{1}+u_{2}(\xi_{1})N_{2}+
a_{11}(\xi_{1}+a_{11}N_{1}+a_{12}N_{2})+a_{12}(\xi_{2}^{\top}+a_{21}N_{1}+a_{22}N_{2})\].
Identifying the tangential and, respectively, normal components on
$M$ we obtain (iv), (v) and (viii) from $(2.6)$. In the same manner,
it result (vi), (vii) and (ix) from $(2.6)$ using $\PP(\PP
N_{2})=N_{2}$.

From $g(PX,Y)=\G(\PP X -u_{1}N_{1}-u_{2}N_{2},Y)=\G(\PP X,Y)$
$=\G(X, \PP Y)=$  $=\G(X,PY+u_{1}(Y)+u_{2}(Y)N_{2})=g(X,PY)$ we get:
the equality (iii) from $(2.7)$. From $\G(\PP X,N_{1})=\G (X, \PP
N_{1})$ we have
\[\G(PX+u_{1}(X)N_{1}+u_{2}(X)N_{2},N_{1})=\G(X,\xi_{1}+a_{11}N_{1}+a_{12}N_{2})\].
Thus, $u_{1}(X)=\G(X,\xi_{1})=g(X,\xi_{1})$ and this yields the
equality (i) from $(2.7)$. In the same manner, using $\G(\PP
X,N_{2})=\G (X, \PP N_{2})$, we obtain (ii) from $(2.7)$.

From $g(PX,Y)=g(X,PY)$, replacing $Y$ with $PY$ we have:
\[g(PX,PY)=g(X,P^{2}Y)=g(X,Y-u_{1}(Y)\xi_{1}-u_{2}(Y)\xi_{2}^{\top}).\]
and from this it results (iv) from $(2.7)$.
\end{proof}

\normalfont From Lemma 1 and Lemma 3 we obtain:

\thm Let $M$ be an n-dimensional submanifold of codimension $1$
isometrically immersed in $\M$, which is also a submanifold of
codimension $1$ and isometrically immersed in an almost product
Riemannian manifold $(\MM,\G,\PP)$. Then, the induced structure on
$M$ by the structure $(\PP,\G)$ from $\MM$ is an $(a, 1)f$
Riemannian structure, determined by $(P, g, u_{1}, u_{2},
\xi_{1}^{\top}, \xi_{2}, (a_{\A\B})_{2})$, (where
$a:=(a_{\A\B})_{2}$ and $f:=P$) which is the same that one induced
on $M$ by the structure $(\overline{P}, \g,u_{1}, \xi_{1}, a_{11})$
(induced on $\M$ by the almost product structure $\PP$ from $\MM$).

\normalfont \medskip We can give a generalization of the Theorem 2.1
as follows:

Let $M:=M_{r}$ be an $n$-dimensional submanifold of codimension $r$
(with $r \geq 2 $) in an almost product Riemannian manifold
$(\MM,\G,\PP)$. We make the following notations: $\MM :=M_{0}$, $\G
:=g^{0}$, $\PP := P_{0}$, such that we have the sequence of
Riemannian immersions given by:
\[(M_{r}, g^{r})\hookrightarrow (M_{r-1}, g^{r-1})
\hookrightarrow ...\hookrightarrow (M_{1}, g^{1})\hookrightarrow
(\MM, \G, \PP)\]  where $g^{i}$ is an induced metric on $M^{i}$ by
the metric $g^{i-1}$ from $M_{i-1}$, ($i \in \{1,...,r\}$) and each
one of $(M_{i},g^{i})$ is a submanifold of codimension $1$,
isometric immersed in the manifold $(M_{i-1},g^{i-1})$ ($i \in
\{1,...,r\}$). Let $i \in \{1,...,r\}$ and $\A_{i}, \B_{i}\in
\{1,...,i\}$). In this condition we obtain:

\thm The $(a, 1)f$ Riemannian structure, determined by the induced
structure $(P_{r}, g^{r}, \xi^{r}_{\A_{r}}, u^{r}_{\A_{r}},
(a^{r}_{\A_{r}\B_{r}}))$
 on an $n$-dimensional submanifold $M:=M_{r}$ of
codimension $r$ (with $r \geq 2 $) in an almost product Riemannian
manifold $(\MM,\G,\PP)$, is the same that one induced on $M$ by any
structures $(P_{i}, g^{i}, \xi^{i}_{\A_{i}}, u_{\A_{i}},
(a^{i}_{\A_{i}\B_{i}}))$ ($i<r$) induced on $M_{i}$ by the almost
product structure $\PP$ on $\MM$, where $f:=P_{r}$ is the tangential
component of $P_{i}$ on $M$, the vector fields $\xi^{r}_{\A_{i}}$ on
$M_{r}$ are the tangential components on $M$ of the tangent vector
fields $\xi^{i}_{\A_{i}}$ from $M_{i}$, the $1$-forms
$u^{r}_{\A_{i}}$ are the restrictions on $M$ of the $1$-forms
$u^{r}_{\A_{i}}$ from $M_{i}$ (for $i<r$), the entries of the $r
\times r$ matrix $a:=(a^{r}_{\A_{r}\B_{r}})$ are defined by
$a^{r}_{\A_{r},\B_{r}}=a^{r}_{\B_{r},\A_{r}}=g^{r}(P_{r-1}(N_{\A_{r}}),N_{\B_{r}})$.

\normalfont

\section{\bf Examples of $(a, 1)f$ Riemannian structures}

\textbf{Example 1.} Let $E^{2p+q}$ be the $(2p+q)$-dimensional
Euclidean space ($p,q \in \mathbb{N}^{*}$). In this example, we
construct an $(\overline{a},1)f$-structure on the sphere
$S^{2p+q-1}(R)\hookrightarrow E^{2p+q}$.

For any point of $E^{2p+q}$ we have its coordinates:
\[(x^{1},...,x^{p},y^{1},...,y^{p},z^{1},...,z^{q}):=(x^{i},y^{i},z^{j})\]
where $i\in \{1,...,p\}$ and $j \in \{1,...,q\}$. The tangent
space $T_{x}(E^{2p+q})$ is isomorphic with $E^{2p+q}$.

Let $\PP:E^{2p+q}\rightarrow E^{2p+q}$ an almost product structure
on $E^{2p+q}$ so that:
\[
\PP(x^{1},...,x^{p},y^{1},...,y^{p},z^{1},...,z^{q})=
(\nu_{1}y^{1},...,\nu_{p}y^{p},\nu_{1}x^{1},...,\nu_{p}x^{p},\varepsilon_{1}z^{1},...,\varepsilon_{q}z^{q})
\leqno(3.1)
\]
 and we use the notation:
\[(\nu_{1}y^{1},...,\nu_{p}y^{p},\nu_{1}x^{1},...,\nu_{p}x^{p},\varepsilon_{1}
z^{1},...,\varepsilon_{q}
z^{q}):=(\nu_{i}y^{i},\nu_{i}x^{i},\varepsilon_{j} z^{j})\] where
$\nu_{i}^{2}=\varepsilon_{j}^{2} = 1$ for every $i\in \{1,...,p\}$
and $j \in \{1,...,q\}$.

The equation of the sphere $S^{2p+q-1}(R)$ is:
\[
\sum_{i=1}^{p}(x^{i})^{2}+\sum_{i=1}^{p}(y^{i})^{2}+
\sum_{j=1}^{q}(z^{j})^{2}=R^{2}
\leqno(3.2)
\]
 where $R$ is its radius and
$(x^{1},...,x^{p},y^{1},...,y^{p},z^{1},...,z^{q}):=(x^{i},y^{i},z^{j})$
are the coordinates of any point of $S^{2p+q-1}(R)$.

We use the following notations:
\[ \sum_{i=1}^{p}(x^{i})^{2}=r_{1}^{2}, \quad \sum_{i=1}^{p}(y^{i})^{2}=r_{2}^{2},
\quad \sum_{j=1}^{q}(z^{j})^{2}=r_{3}^{2}\] and
$r_{1}^{2}+r_{2}^{2}=r^{2}$. Thus we have $r^{2}+r_{3}^{2}=R^{2}$.

We remark that an unit normal vector field on sphere
$S^{2p+q-1}(R)$ has the form:
\[
N_{1}:= \frac{1}{R}(x^{i},y^{i},z^{j}), \leqno(3.3)
\]
for $i \in \{1,...,p\}$ and $j \in \{1,...,q\}$ and we have $\PP
 N_{1}=\frac{1}{R}(\nu_{i}y^{i}, \nu_{i}x^{i},
 \varepsilon_{j}z^{j})$.

For any tangent vector field:
\[\overline{X}=(X^{1},...,X^{p},Y^{1},...,Y^{p},Z^{1},...,Z^{q}):=(X^{i},Y^{i},Z^{j})\]
on $S^{2p+q-1}(R)$ we have:
\[
\sum_{i=1}^{p}x^{i}X^{i}+\sum_{i=1}^{p}y^{i}Y^{i}+\sum_{j=1}^{q}z^{j}Z^{j}=0, \leqno(3.4)
\]

From $(1.1)$ and $(1.2)$ we have the decompositions of $\PP
\overline{X}$ and $\PP N_{1}$ in tangential and normal components,
respectively, at the sphere $S^{2p+q-1}(R)$.

In the following issue we use the notations $\overline{a}:=a_{11}$
and $f:=\overline{P}$:
\[
 \sigma = \sum_{i=1}^{p}\nu_{i}x^{i}y^{i},
\quad \tau=\sum_{j=1}^{q}\varepsilon_{j}(z^{j})^{2}, \leqno(3.5)
\]

\[
\gamma = \sum_{i=1}^{p}\nu_{i}(x^{i}Y^{i}+y^{i}X^{i}), \quad
\mu=\sum_{j=1}^{q}\varepsilon_{j}z^{j}Z^{j}
\leqno(3.6)
\]
for any point $(x^{i}, y^{i}, z^{j})$ of $S^{2p+q-1}(R)$ and for any
tangent vector field $\overline{X}=(X^{i}, Y^{i}, Z^{j})$ ($i\in
\{1,...,p\}$ and $j \in \{1,...,q\}$). Using the first Lemma, we
obtain an $(\overline{a}, 1)f$ structure on the sphere
$S^{2p+q-1}(R)\hookrightarrow E^{2p+q}$ (with $\overline{g}:=<>$),
determined by $(\PB, <>, \xi_{1}, u_{1}, a_{11})$ which has the
elements as follows:
\[
a_{11}=\frac{2 \sigma + \tau}{R^{2}},\leqno(3.7)
\]
\[
u_{1}(\overline{X})=\gamma+\tau, \leqno(3.8)
\]
\[
\xi_{1}=\frac{1}{R}(\nu_{i}y^{i}-a_{11}x^{i},
\nu_{i}x^{i}-a_{11}y^{i}, (\varepsilon_{j}-a_{11})z^{j}),
\leqno(3.9)
\]
and:
\[
\overline{P}(\overline{X})=(\nu_{i}Y^{i}-\frac{u_{1}(\overline{X})}{R}x^{i},
\nu_{i}X^{i}- \frac{u_{1}(\overline{X})}{R}y^{i},
\varepsilon_{j}Z^{j}-\frac{u_{1}(\overline{X})}{R}z^{j}).
\leqno(3.10)
\]

\smallskip

\textbf{Example 2.} In this example, we construct an $(a,
1)f$-structure on the product of spheres $S^{2p-1}(r) \times
S^{q-1}(r_{3})$. Let $E^{2p+q}$ ($p,q \in \mathbb{N}^{*}$) be the
Euclidean space ($p,q \in \mathbb{N}^{*}$) endowed with the almost
product Riemannian structure $\PP$ defined in $(3.1)$. We set
$E^{2p+q}=E^{2p} \times E^{q}$ and in each of spaces $E^{2p}$ and
$E^{q}$ respectively, we consider the spheres:
\[S^{2p-1}(r)=\{(x^{1},...,x^{p},y^{1},...,y^{p}), \sum_{i=1}^{p}((x^{i})^{2}+(y^{i})^{2})=r^{2}\}\]
and respectively:
\[ S^{q-1}(r_{3})=\{(z^{1},...,z^{q}), \sum_{j=1}^{q}(z^{j})^{2}=r_{3}^{2}\}\]
 where $r^{2}+r_{3}^{2}=R^{2}$. Any point of the
product manifold $S^{2p-1}(r) \times S^{q-1}(r_{3})$ has the
coordinates
$(x^{1},...,x^{p},y^{1},...,y^{p},z^{1},...,z^{q}):=(x^{i},y^{i},z^{j})$
which verify $(3.2)$. Thus $S^{2p-1}(r) \times S^{q-1}(r_{3})$ is a
submanifold of codimension $2$ in $E^{2p+q}$. Furthermore,
$S^{2p-1}(r) \times S^{q-1}(r_{3})$ is a submanifold of codimension
$1$ in $S^{2p+q-1}(R)$. Therefore, we have:
\[S^{2p-1}(r) \times S^{q-1}(r_{3})
 \hookrightarrow S^{2p+q-1}(R) \hookrightarrow E^{2p+q} \].
The tangent space in a point $(x^{i}, y^{i}, z^{j})$ at the product
of spheres $S^{2p-1}(r) \times S^{q-1}(r_{3})$ is
$T_{(x^{1},...,x^{p},y^{1},...,y^{p},\underbrace{o,...,o}_{q})}
S^{2p-1}(r)
 \oplus
 T_{(\underbrace{o,...,o}_{2p},z^{1},...,z^{q})}S^{q-1}(r_{3})$.

A vector $(X^{1},...,X^{p},Y^{1},...,Y^{p})$ from
$T_{(x^{1},...,x^{p},y^{1},...,y^{p})}E^{2p}$ is tangent to
$S^{2p-1}(r)$ if and only if:
\[
\sum_{i=1}^{p}x^{i}X^{i}+\sum_{i=1}^{p}y^{i}Y^{i}=0 \leqno(3.11)
\]
and it can be identified with
$(X^{1},...,X^{p},Y^{1},...,Y^{p},\underbrace{0,...,0}_{q})$ from
$E^{2p+q}$. A vector $(Z^{1},...,Z^{q})$ from
$T_{(z^{1},...,z^{q})}E^{q}$ is tangent to $S^{q-1}(r_{3})$ if and
only if:
\[
\sum_{j=1}^{q}z^{j}Z^{j}=0 \leqno(3.12)
\] and it can be identified with
$(\underbrace{0,...,0}_{2p},Z^{1},...,Z^{q})$ from $E^{2p+q}$.

Consequently, for any point $(x^{i},y^{i},z^{j}) \in S^{2p-1}(r)
\times S^{q-1}(r_{3})$ we have $(X^{i},Y^{i},Z^{j}) \in
T_{(x^{1},...,x^{p},y^{1},...,y^{p},z^{1},...,z^{q})}(S^{2p-1}(r)
\times S^{q-1}(r_{3}))$ if and only if the equations $(3.11)$ and
$(3.12)$ are satisfied. Furthermore, we remark that
$(X^{i},Y^{i},Z^{j})$ is a tangent vector field at $S^{2p+q-1}(R)$
and from this it follows that:
\[T_{(x^{i},y^{i},z^{j})}(S^{2p-1}(r) \times S^{q-1}(r_{3})) \subset
T_{(x^{i},y^{i},z^{j})}S^{2p+q-1}(R),\] for any point
$(x^{i},y^{i},z^{j}) \in S^{2p-1}(r) \times S^{q-1}(r_{3})$.

The normal unit vector field $N_{1}$  at $S^{2p+q-1}(R)$ given by
$(3.3)$ is also a normal vector field at $(S^{2p-1}(r) \times
S^{q-1}(r_{3}))$ when it is considered in its points. We construct
an unit vector field $N_{2}$ on $S^{2p+q-1}$ as follows:
\[
N_{2}=\frac{1}{R}(\frac{r_{3}}{r}x^{i},\frac{r_{3}}{r}y^{i},-\frac{r}{r_{3}}z^{j})
\leqno(3.13)
\]

It is obvious that $(N_{1},N_{2})$ defined in $(3.3)$ and $(3.13)$
is a local orthonormal basis in
$T_{(x^{i},y^{i},z^{j})}^{\bot}S^{2p-1}(r) \times S^{q-1}(r_{3})$ in
any point $(x^{i},y^{i},z^{j}) \in S^{2p-1}(r) \times
S^{q-1}(r_{3})$. Using Lemma 2 and Lemma 3, we obtain the structure
$(\widehat{P}, <>, \widehat{\xi}_{1},
\widehat{\xi}_{2},\widehat{u}_{1}, \widehat{u}_{2}, \widehat{a})$ on
the product of spheres $S^{2p-1}(r) \times S^{q-1}(r_{3})$, induced
by the almost product Riemannian structure $(\PP, <,>)$ as
follows:\\
$\cdot$ the matrix $a:=(a_{\A\B})_{2}$ is given by:
\[
a:= \begin{pmatrix}
   \frac{2 \sigma + \varepsilon r_{3}^{2}}{R^{2}} &  \frac{(2 \sigma - \varepsilon r^{2})r_{3}}{rR^{2}}\\
   \frac{(2 \sigma - \varepsilon r^{2})r_{3}}{rR^{2}} & \frac{2 \sigma r_{3}^{2}+ \varepsilon r^{4}}{r^{2}R^{2}}
\end{pmatrix},
\leqno(3.14)
\]
$\cdot$ the tangent vector fields have the form:
\[
\xi_{1}=\frac{1}{R} (\nu_{i}y^{i}-\frac{2 \sigma}{r^{2}}x^{i},
\nu_{i}x^{i}-\frac{2
\sigma}{r^{2}}y^{i},(\varepsilon_{j}-\frac{\tau}{r_{3}^{2}})z^{j}),
\leqno(3.15)
\]
and:
\[
\xi_{2}=\frac{1}{R} (\frac{r_{3}}{r}(\nu_{i}y^{i}-\frac{2 \sigma
}{r^{2}}x^{i}), \frac{r_{3}}{r}(\nu_{i}x^{i}-\frac{2
\sigma}{r^{2}}y^{i}),-\frac{r}{r_{3}}((\varepsilon_{j}-\frac{\tau}{r_{3}^{2}})z^{j}),
\leqno(3.16)
\]
$\cdot$ the $1$-forms are given by:
\[
u_{1}(X)=\frac{1}{R}(\gamma + \mu), \quad
u_{2}(X)=\frac{1}{R}(\frac{r_{3}}{r}\gamma - \frac{r}{r_{3}} \mu),
\leqno(3.17)
\]
and the $(1,1)$ tensor field P has the form:
\[
 P(X)=(\nu_{i}Y^{i}-\frac{\gamma}{r^{2}}x^{i},\nu_{i}X^{i}-\frac{\gamma}{r^{2}}y^{i},\varepsilon_{j}Z^{j}-\frac{\mu}{r_{3}^{2}}z^{j})
\leqno(3.18)
\]
for any tangent vector field
 $X:=(X^{i},Y^{i},Z^{j}) \in
T_{(x^{1},...,x^{p},y^{1},...,y^{p},z^{1},...,z^{q})}(S^{2p-1}(r)
\times S^{q-1}(r_{3}))$ and any point $(x^{i},y^{i},z^{j}) \in
S^{2p-1}(r) \times S^{q-1}(r_{3})$. For $a:=(a_{\A\B})_{2}$ and
$f:=P$, the structure $(\widehat{P}, <>, \widehat{\xi}_{1},
\widehat{\xi}_{2},\widehat{u}_{1}, \widehat{u}_{2}, \widehat{a})$ is
an $(a, 1)f$ Riemannian structure induced on the on the product of
spheres $S^{2p-1}(r) \times S^{q-1}(r_{3})$ which is a submanifold
of codimension 2 in the Euclidean space $E^{2p+q}$.

\normalfont

\end{document}